\documentclass[12pt]{amsart}

\usepackage{graphicx,amssymb,latexsym,amsfonts,txfonts,amsmath,amsthm}
\usepackage{pdfsync,color,tabularx,rotating}
\usepackage[all,cmtip]{xy}
\usepackage{url}
\usepackage{tikz}
\usepackage{amssymb}
\usepackage{comment}

\advance\hoffset by -0.9 truecm

\newtheorem{theo}{Theorem}[section]

\newtheorem{conj}[theo]{Conjecture}
\newtheorem{prob}[theo]{Problem}

\theoremstyle{remark}

\theoremstyle{remark}

\theoremstyle{remark}
\newtheorem{defi}[theo]{\bf Definition}

\theoremstyle{remark}

\begin{document}

\title{A number-theoretic problem concerning pseudo-real Riemann surfaces}

\author{Gareth A. Jones and Alexander K. Zvonkin}

\address{School of Mathematical Sciences, University of Southampton, Southampton SO17 1BJ, UK}
\email{G.A.Jones@maths.soton.ac.uk}

\address{LaBRI, Universit\'e de Bordeaux, 351 Cours de la Lib\'eration, F-33405 Talence Cedex, France}

\email{zvonkin@labri.fr}

\subjclass[2010]{11A41, 11M06, 20B25, 30F10} 

\keywords{Riemann surface, pseudo-real, prime number, Riemann zeta function, Dirichlet series, Bateman--Horn Conjecture}

\begin{abstract}
Motivated by their research on automorphism groups of pseudo-real Riemann surfaces, Bujalance, Cirre and Conder have conjectured that there are infinitely many primes $p$ such that $p+2$ has all its prime factors $q\equiv -1$ mod~$(4)$. We use theorems of Landau and Raikov to prove that the number of integers $n\le x$ with only such prime factors $q$ is asymptotic to $cx/\sqrt{\ln x}$ for a specific constant $c=0.4865\ldots$. Heuristic arguments, following Hardy and Littlewood, then yield a conjecture that the number of such primes $p\le x$ is asymptotic to $c'\int_2^x(\ln t)^{-3/2}dt$ for a constant $c'=0.8981\ldots$. The theorem, the conjecture and a similar conjecture applying the Bateman--Horn Conjecture to other pseudo-real Riemann surfaces are supported by evidence from extensive computer searches.
\end{abstract}

\maketitle


\section{Introduction}\label{sec:Intro}

A compact Riemann surface is real (meaning definable, as a complex algebraic curve, over $\mathbb R$) if and only if it has an orientation-reversing automorphism of order $2$; it is called {\sl pseudo-real\/} if it has an orientation-reversing automorphism, but none of order $2$.
In~\cite{BCC}, Bujalance, Cirre and Conder have proved that for each $g\ge 2$, the maximum order $M^+_{\rm ab}(g)$ of an abelian group of (orientation-preserving) automorphisms of a pseudo-real Riemann surface of genus $g$ is at least $g$. To demonstrate the sharpness of this lower bound, their Theorem 4.8 presents a set  $\mathcal A$ of integers $g$ for which the cyclic group ${\rm C}_g$ of order $g$ attains the upper bound $M^+_{\rm ab}(g)=g$. These have the form $g=p+1$ where $p$ is
what we will call a BCC prime, defined as follows:
\begin{defi}
A BCC prime is a prime number $p$ such that $p+2$ has only prime factors $q\equiv -1$ mod~$(4)$.
\end{defi}
\noindent The first sixteen BCC primes are
\begin{equation}\label{eq:smallBCC}
5, 7, 17, 19, 29, 31, 41, 47, 61, 67, 79, 97, 101, 127, 131\;\hbox{and}\; 137.
\end{equation}

The authors of~\cite{BCC} show that ${\rm C}_g$ is the orientation-preserving subgroup of ${\rm C}_{2g}$ in its action on a compact Riemann surface $\mathcal S$ of genus $g$, with ${\mathcal S}/{\rm C}_{2g}$ having signature $(1; - ; [2, 2, g]; \{-\})$. (This is both a description of the quotient space as an orbifold, namely the real projective plane with cone points of orders $2, 2$ and $g$, and an encoded presentation of the non-euclidean crystallographic group which uniformises it.)
If $g\in{\mathcal A}$ then $g$ is even, so the only involution in ${\rm Aut}\,{\mathcal S}$ preserves orientation, and hence $\mathcal S$ is pseudo-real. They then use a detailed case-by-case argument to show that if also $g > 30$ then no pseudo-real Riemann surface of genus $g$ has a larger abelian group of automorphisms. For example, since $p+2$ has no prime factors $q\equiv 1$ mod~$(4)$, ${\rm Aut}\,{\rm C}_{p+2}$ (of order $\phi(p+2)$) has no elements of order $4$, a fact used in subcase (2c) of the proof of Theorem~4.8 to eliminate ${\rm C}_{g+1}={\rm C}_{p+2}$. Although our main emphasis here is on number theory and computation, for readers interested in Riemann surfaces we have described the construction of the surfaces $\mathcal S$ in more detail in an appendix, Section~\ref{sec:pseudo-real}.

In~\cite{BCC}, $\mathcal A$ is described as `a very large and possibly infinite set'.
If $\mathcal A$ is finite, then conceivably there is a better lower bound for $M^+_{\rm ab}(g)$, valid for all sufficiently large $g$ (in particular, larger than all those in $\mathcal A$, which can then be regarded as small exceptions). To avoid this possibility, it is important to know whether $\mathcal A$ is infinite.
We are therefore interested in the following problem:

\begin{prob}\label{pro:BCC}
Are there infinitely many BCC primes, or equivalently, is the set $\mathcal A$ infinite?
\end{prob}
\noindent Of course, the answer will be positive if one can prove that there are infinitely many prime pairs $p, p+2$ with $p\equiv 1$ mod~$(4)$, but even without this extra congruence condition, the existence or otherwise of infinitely many prime pairs is a very difficult open problem. Nevertheless, this would follow from a proof of various conjectures in Number Theory, such as the Bateman--Horn Conjecture~\cite{BH} or Schinzel's Hypothesis H~\cite{SS}. Here we give a heuristic argument, supported by computational evidence for their asymptotic distribution, for the following

\begin{conj}\label{con:BCC}
There are infinitely many BCC primes.
\end{conj}

(It follows from a theorem of Iwaniec~\cite{Iwa} that there are infinitely many primes $p$ such that $p+2$ has only prime factors $q\equiv 1$ mod~$(4)$, whereas the corresponding result with $q\equiv -1$ mod~$(4)$ is unproved; see the Digression in Section~\ref{sec:Distribution} for details.)

In support of Conjecture~\ref{con:BCC} we use theorems of Landau~\cite{Lan} and of Raikov~\cite{Rai} to give, in Theorem~\ref{th:g-dist}, an asymptotic estimate of the form $cx/\sqrt{\ln x}$ (with $c$ a specific constant) for the number of integers $n\le x$ with only prime factors $q\equiv -1$ mod~$(4)$. Heuristic arguments, inspired by Hardy and Littlewood~\cite{HL}, then allow us to make Conjecture~\ref{conj:a-distr}, that the number of BCC primes $p\le x$ is asymptotic to
\[c'\negthickspace\int_2^x\negthickspace\frac{dt}{(\ln t)^{3/2}}\]
for a specific constant $c'$. We give numerical evidence, based on computer searches, to support both Theorem~\ref{th:g-dist} and Conjecture~\ref{conj:a-distr}.
In addition, in Section~\ref{sec:Similar}, we briefly consider another set of primes, defined by similar but more complicated conditions, also arising from a construction in~\cite{BCC}.

The formulae we obtain involve various multiplicative factors. These are defined in terms of well-known constants such as $e$, $\gamma$ and $\pi$, together with constants such as $C(k,u)$ (defined later) which have been computed elsewhere to over $100$ decimal places, so in principle they can be computed with similar accuracy. However, for simplicity we have in most cases presented numerical data to about ten significant figures, since this is adequate for the arguments we wish to present. 


\section{Raikov's Theorem}\label{sec:Raikov}

If $\mathcal P$ is a set of prime numbers, then let us define an integer $n\in\mathbb N$ to be a $\mathcal P$-{\sl integer\/} if all its prime factors are elements of $\mathcal P$, and let ${\mathcal P}^*$ denote the set of all $\mathcal P$-integers. If we define $g_n=1$ or $0$ as $n\in{\mathcal P}^*$ or not, then the function
\begin{equation}
g(x):=\sum_{n\le x}g_n
\end{equation}
gives the number of $\mathcal P$-integers $n\le x$. We will consider the asymptotic behaviour of $g(x)$ as $x\to\infty$.

In~\cite{Lan} Landau showed that if, for some $k\in{\mathbb N}$, $\mathcal P$ is the set of all primes in the union of $l$ distinct congruence classes of units mod~$(k)$, then
\[g(x)\sim \frac{cx}{(\ln x)^{1-l/\phi(k)}}\quad\hbox{as}\quad x\to\infty\]
for some constant $c>0$. In the context of Problem~\ref{pro:BCC}, if we take $k=4$ and define ${\mathcal P}$ to be the set of primes in the congruence class $[-1]\in{\mathbb Z}^*_4$, we see that
\[g(x)\sim \frac{cx}{\sqrt{\ln x}}\quad\hbox{as}\quad x\to\infty\]
for some $c>0$. Our aim is to determine this constant $c$, and then to compare the resulting estimates with the actual values of $g(x)$ for various $x$.

For any set $\mathcal P$ of primes, the function $n\mapsto g_n$ is completely multiplicative, so the corresponding Dirichlet series
\begin{equation}
F(s):=\sum_{n=1}^{\infty}\frac{g_n}{n^s}
\end{equation}
has an Euler product expansion
\begin{equation}
F(s)=\prod_{q\in{\mathcal P}}\left(1-\frac{1}{q^s}\right)^{-1}.
\end{equation}
By comparison with the Riemann zeta function
\begin{equation}
\zeta(s):=\sum_{n=1}^{\infty}\frac{1}{n^s}=\prod_{q\;{\rm prime}}\left(1-\frac{1}{q^s}\right)^{-1},
\end{equation}
this series and product for $F(s)$ converge for all $s\in\mathbb C$ with ${\rm Re}(s)>1$.

A theorem of Raikov~\cite{Rai}, as given in~\cite[Theorem~2.4.1]{CoMu}, states the following:
\begin{theo}
Let $F(s)=\sum_{n\ge 1}a_n/n^s$ be a Dirichlet series with non-negative coefficients, converging for ${\rm Re}(s) > 1$. 
Suppose that $F(s)$ extends analytically at all points on ${\rm Re}(s) = 1$ apart from $s=1$, and that at $s=1$ we can write
\begin{equation}\label{eq:H-defn}
F(s)=\frac{H(s)}{(s-1)^{1-\alpha}}
\end{equation}
for some $\alpha\in{\mathbb R}$ and some $H(s)$ holomorphic in the region ${\rm Re}(s)\ge 1$ and nonzero there. Then
\begin{equation}\label{eq:g-est}
\sum_{n\le x}a_n\sim\frac{cx}{(\ln x)^{\alpha}}
\end{equation}
as $x\to\infty$, with
\begin{equation}
c=\frac{H(1)}{\Gamma(1-\alpha)}
\end{equation}
where $\Gamma$ is the Gamma function
\[\Gamma(z)=\int _0^{\infty }\negthickspace t^{z-1}e^{-t}\,dt\]
for ${\rm Re}(z)>0$.
\end{theo}

In simple terms, (\ref{eq:H-defn}) means that $F(s)$ has a pole of order $1-\alpha$ at $s=1$, in the sense that it is a non-zero holomorphic multiple of $1/(s-1)^{1-\alpha}$ near $s=1$, and $H(1)$ is the limiting value of the multiplier as $s\to 1$.

In the next section we will determine the function $H$ and hence the constants $\alpha$ and $c$ for the sequence $(g_n)$.
We will do this by expressing the corresponding function $F$ in terms of $\zeta$ and related functions, so that the required analyticity conditions for $F$ will follow from those of $\zeta$.


\section{Distribution of $\mathcal P$-integers\label{sec:Distribution}}

In the case of Problem~\ref{pro:BCC} we take
\[{\mathcal P}:=\{q\mid q\;\hbox{is prime and}\; q\equiv -1\;{\rm mod}\,(4)\}.\]
Thus the $\mathcal P$-integers $n\le 100$ are
\[1, 3, 7, 9, 11, 19, 21, 23, 27, 31, 33, 43, 47,
49, 57, 59, 63, 67, 69, 71, 77, 79, 81, 83, 93, 99.\]
Note that we include $1$ as the empty product of primes $q\equiv -1$ mod~$(4)$. These $26$ integers represent just over half of the odd integers $n\le 100$. However, this proportion decreases towards $0$ as the upper bound increases: see Theorem~\ref{th:g-dist} and Table~\ref{tab:g(x)est}, later in this section.

In order to determine the corresponding constants $c$ and $\alpha$ we can determine the right-hand side of (\ref{eq:H-defn}) by comparing $F(s)$, or more precisely $F(s)^2$, with $\zeta(s)$, which is meromorphic with a simple pole of residue $1$ at $s=1$ (see~\cite[Theorem 8.28]{EW}, for example).

For ${\rm Re}(s)>1$ we have
\[\frac{F(s)^2}{\zeta(s)}=\left(1-\frac{1}{2^s}\right)\cdot\prod_{q\equiv -1(4)}\left(1-\frac{1}{q^s}\right)^{-2}\cdot\prod_{q\equiv\pm 1(4)}\left(1-\frac{1}{q^s}\right)
=\left(1-\frac{1}{2^s}\right)\cdot\prod_{q\equiv \pm 1(4)}\left(1-\frac{1}{q^s}\right)^{\pm 1}.\]
so that
\begin{equation}\label{eq:F-eqn}
(s-1)F(s)^2=\left(1-\frac{1}{2^s}\right)\cdot\prod_{q\equiv \pm 1(4)}\left(1-\frac{1}{q^s}\right)^{\pm 1}\cdot(s-1)\zeta(s).
\end{equation}
As $s\to 1$ we have
\[\left(1-\frac{1}{2^s}\right)\to\frac{1}{2},\quad
\prod_{q\equiv \pm 1(4)}\left(1-\frac{1}{q^s}\right)^{\pm 1}\to\prod_{q\equiv \pm 1(4)}\left(1-\frac{1}{q}\right)^{\pm 1}
\quad\hbox{and}\quad (s-1)\zeta(s)\to 1.\]

In order to deal with the second of these limits, we can use some results of Languasco and Zaccagnini~\cite{LZ1, LZ2, LZ3}.
For each integer $k\ge 3$ and each unit $u$ mod~$(k)$ they define a non-zero Mertens-type constant $C(k,u)$ by the asymptotic estimate
\begin{equation}\label{eq:P-defn}
P(k,u;x):=\prod_{x\ge q\equiv u(k)}\left(1-\frac{1}{q}\right)=\frac{C(k,u)}{(\ln x)^{1/\phi(k)}}+O\left(\frac{1}{(\ln x)^{1+1/\phi(k)}}\right)\quad\hbox{as}\quad x\to\infty,
\end{equation}
where the product is over all primes $q\equiv u$ mod~$(k)$ such that $q\le x$,
and $\phi$ is Euler's totient function.
In~\cite[Equation~(2)]{LZ2} they show that
\begin{equation}\label{eq:C-eqn}
C(k,u)^{\phi(k)}=e^{-\gamma}\prod_{q\,\rm prime}\left(1-\frac{1}{q}\right)^{\alpha(q)}
\end{equation}
where $\gamma$ is the Euler--Mascheroni constant, the product is now over all primes $q$, and $\alpha(q):=\phi(k)-1$ or $-1$ as $q\equiv u$ mod~$(k)$ or not.
By taking $k=4$ and $u=1$ in (\ref{eq:C-eqn}) we see that
\[\prod_{q\equiv \pm 1(4)}\left(1-\frac{1}{q}\right)^{\pm 1}=e^{\gamma}C^2\left(1-\frac{1}{2}\right)=\frac{e^{\gamma}C^2}{2}\]
where $C:=C(4,1)$.

 It therefore follows by taking limits as $s\to 1$ in~(\ref{eq:F-eqn}) that
\[H(1)^2=\lim_{s\to 1}\,(s-1)F(s)^2=\frac{e^{\gamma}}{4}C^2\ne 0,\]
so that this particular instance of (\ref{eq:H-defn}) becomes $F(s)=H(s)/(s-1)^{1/2}$ and hence
\[\alpha=\frac{1}{2}\quad\hbox{and}\quad H(1)=\sqrt{e^{\gamma}}\cdot\frac{C}{2}.\]
Since $\Gamma(\frac{1}{2})=\sqrt\pi$ we therefore have
\begin{equation}\label{eq:c-formula}
c=\frac{H(1)}{\Gamma(\frac{1}{2})}=\sqrt{\frac{e^{\gamma}}{\pi}}\cdot\frac{C}{2},
\end{equation}\label{eq:g(x)-estimate}
giving the following result\footnote{Here we cannot resist expressing our pleasure at the appearance of the constant $\sqrt{e^{\gamma}/\pi}$, in which the three basic constants $e$, $\pi$ and $\gamma$ of Analysis are united by three of the basic operations (division, exponentiation and taking square roots) of Algebra; its application to Number Theory adds to the pleasure.}:

\begin{theo}\label{th:g-dist}
The function $g(x)$ satisfies

\begin{equation}\label{eq:g-estimate}
g(x)\sim\frac{cx}{\sqrt{\ln x}}\;=\sqrt{\frac{e^{\gamma}}{\pi}}\cdot\frac{C}{2}\cdot\frac{x}{\sqrt{\ln x}}\quad\hbox{as}\quad x\to\infty.
\end{equation}
\end{theo}

 In~\cite{LZ2, LZ3} Languasco and Zaccagnini have evaluated many of these constants $C(k,u)$ to over $100$ decimal places. These include
\begin{equation}\label{eq:C-value}
C=C(4,1)=1.2923041571286886071091383898704320653429\ldots,
\end{equation}
which, together with
\[\sqrt{\frac{e^{\gamma}}{\pi}}=0.7529495060464205959354997575876048108386\ldots,\]
allows us to evaluate
\begin{equation}\label{eq:c_0-value}
c=0.4865198883858909971272456405868234055382\ldots
\end{equation}

We use this in Table~\ref{tab:g(x)est} to test the accuracy of the resulting estimates given by Theorem \ref{th:g-dist}. The second column gives the number $g(x)$ of ${\mathcal P}$-integers $n\le x$ for $x=10^k$ with $k=1,\ldots, 9$. For example, the $25$ such integers $n\le 10 ^2$ are
\[3,\, 7,\, 9,\, 11,\, 19,\, 21,\, 23,\, 27,\, 31,\, 33,\, 43,\, 47,\, 49,\, 57,\, 59,\, 63,\, 67,\, 69,\, 71,\, 77,\, 79,\, 81,\, 83,\, 93,\, 99.\]
Just $12$ of these correspond to BCC primes $p=n-2$, namely those primes displayed in Section~\ref{sec:Intro},
with $12$ associated genera $g=p+1=n-1\in{\mathcal A}$. The third column of Table~\ref{tab:g(x)est} gives the estimates $cx/\sqrt{\ln x}$ for $g(x)$, rounded to the nearest integer, and the fourth column gives their errors.

\begin{table}[htbp]
\begin{center}
\begin{tabular}{c|c|c|c|c|c}
$x$ 	& $g(x)$ & $cx/\sqrt{\ln x}$ & error & $E(x)$ & error \\
\hline
10 	     & 4 	           & 3.21 	        & $-19.75\,\%$             & 3.08 	            & $-23.00\,\%$	\\
$10^2$  & 26	           & 22.67 	        & $-12.81\,\%$        & 25.58 	   & $-1.62\,\%$	\\
$10^3$  & 201	           & 185.11        & $-7.91\,\%$        & 202.61 	   & 1.80\,\%	\\
$10^4$  & 1\,680         & 1\,603.11 	& $-4.58\,\%$        & 1\,710.35     & 1.80\,\% 	\\
$10^5$  & 14\,902         & 14\,338.63 	& $-3.78\,\%$        & 15\,069.0     & 1.12\,\%	\\
$10^6$  & 135\,124 	   & 130\,893.21                  & $-3.13\,\%$   & 136\,274.75 	      & 0.85\,\%	\\
$10^7$  & 1\,243\,370         & 1\,211\,835.68        & $-2.54\,\%$    & 1\,253\,639.87   & 0.83\,\%	\\
$10^8$  & 11\,587\,149 	& 11\,335\,684.78 	 & $-2.17\,\%$    & 11\,672\,710.45  & 0.74\,\% \\
$10^9$  & 108\,941\,388 	& 106\,873\,861.02 	 & $-1.90\,\%$    & 109\,666\,579.94	 & 0.67\,\% \\
$10^{10}$  & 1\,031\,330\,156  & 1\,013\,894\,469.43 & $-1.69\,\%$ & 1\,037\,530\,754.16 & 0.60\,\%
\end{tabular}
\medskip
\caption{Estimates $cx/\sqrt{\ln x}$ and $E(x)=c\int_2^x\frac{dt}{\sqrt{\ln t}}$ for $g(x)$.}
\label{tab:g(x)est}
\end{center}
\end{table}

An error of $1.69\%$ for $x=10^{10}$ is not impressive, but then the simple form $\pi(x)\sim x/\ln x$ of the Prime Number Theorem, which can also be viewed as an instance of Raikov's Theorem, has an error of about $5\%$ here. Motivated by the much greater accuracy of ${\rm Li}(x)=\int_2^x(\ln t)^{-1}dt$ as an estimate for $\pi(x)$, we considered whether
\[E(x):=c\negthickspace\int_2^x\negthickspace\frac{dt}{\sqrt{\ln t}}\]
might give a better estimate for $g(x)$. The evidence is given in the fifth and sixth columns of Table~\ref{tab:g(x)est}.

More generally, Raikov's Theorem (\ref{eq:g-est}) can be restated as
\[\sum_{n\le x}a_n\sim c\negthickspace\int_2^x\negthickspace\frac{dt}{(\ln t)^{\alpha}},\]
where $c$ and $\alpha$ are as defined in Section~\ref{sec:Raikov},
since the two estimates are asymptotically equivalent. 
Our experience, together with heuristic arguments based on the expected values of certain random variables,
suggests that this integral form will usually give more accurate approximations.

\medskip

{\bf Digression.} Motivated by curiosity rather than problems involving Riemann surfaces, we note that if one defines
\[{\mathcal P}^+:=\{q\mid q\;\hbox{is prime and}\; q\equiv 1\;{\rm mod}\,(4)\}\]
then by essentially the same argument the corresponding counting function $g^+(x)$ for $\mathcal P^+$-integers $n\le x$ is estimated by a similar formula
\[g^+(x)\sim\frac{c^+x}{\sqrt{\ln x}}\;=\sqrt{\frac{e^{\gamma}}{\pi}}\cdot\frac{C^+}{2}\cdot\frac{x}{\sqrt{\ln x}}\quad\hbox{as}\quad x\to\infty,\]
where
\[C^+:=C(4,-1)=0.8689277682343238299091527791046529122939\ldots\]
according to~\cite{LZ2,LZ3}. Comparing this value with that of $C$ shows that $\mathcal P$-integers appear nearly $50\%$ more frequently than $\mathcal P^+$-integers. This is despite the fact that the numbers $\pi^{\pm}(x)$ of primes $q\le x$ satisfying $q\equiv\pm 1$ mod~$(4)$ satisfy
\[\pi^+(x)\sim\pi^-(x)\sim\frac{1}{2}\pi(x)\quad\hbox{as}\quad x\to\infty.\]
The reason for this surprisingly large bias is that the smallest primes $q\equiv 1$ mod~$(4)$, such as $5, 13, 17, 29,\ldots$, are almost always larger than the corresponding primes congruent to $-1$ mod~$(4)$, such as  $3, 7, 11, 19,\ldots$ (see~\cite{GM} for a very readable account of this and other similar phenomena), so that $\mathcal P^+$-integers tend to be larger and hence less frequently found than $\mathcal P$-integers.
For example, the first ten $\mathcal P$-integers are $3, 7, 9, 11, 19, 21, 23, 27, 31, 33,$
while the corresponding $\mathcal P^+$-integers are $5, 13, 17, 25, 29, 37, 41, 53, 61, 65$.

Table~\ref{tab:good-vs-bad} gives the numbers of $\mathcal P$- and $\mathcal P^+$-integers $n\le x$ for various $x$.
The constants $C(4,1)=1.29230415713\ldots$ and $C(4,-1)=0.86892776823\ldots$ of interest to us have ratio $C(4,-1)/C(4,1)=0.6723864220\ldots$.
It looks plausible that the ratios in the last column of the table tend to this value.

\begin{table}[htbp]
\begin{center}
\begin{tabular}{c|c|c|c}
$x$		& $g(x)$	& $g^+(x)$ 	& $g^+(x)/g(x)$ \\
\hline
10		& 4				& 2 				& 0.5000000000	\\
$10^2$	& 26				& 15				& 0.5769230769 	\\
$10^3$	& 201			& 123			& 0.6119402985	\\
$10^4$	& 1\,680			& 1\,074			& 0.6392857142	\\
$10^5$	& 14\,902  		& 9\,623			& 0.6457522480	\\
$10^6$	& 135\,124		& 87\,882		        & 0.6503803913	\\
$10^7$	& 1\,243\,370    	& 814\,183		& 0.6548195629	\\
$10^8$	& 11\,587\,149	        & 7\,618\,317	        & 0.6574798511	\\
$10^9$	& 108\,941\,388	& 71\,838\,469	        & 0.6594231110	\\
\end{tabular}
\medskip
\caption{Numbers $g(x)$ and $g^+(x)$ of $\mathcal P$- and $\mathcal P^+$-integers $n\le x$.}
\label{tab:good-vs-bad}
\end{center}
\end{table}

It follows from a theorem of Iwaniec~\cite{Iwa} that for each integer $A\ne 0$ the number of primes $p\le x$ of the form $\xi^2+\eta^2+A$ with $\gcd(\xi,\eta)=1$ has order of magnitude $x/(\log x)^{3/2}$. Now a positive integer $n$ is properly represented by the quadratic form $\xi^2+\eta^2$, that is, with $\gcd(\xi,\eta)=1$, if and only if $n$ is not divisible by $4$ or by any prime $q\equiv-1$ mod~$(4)$, or equivalently $n=m$ or $2m$ where $m$ is a ${\mathcal P}^+$-integer. Taking $A=-2$ we see that there are infinitely many primes $p$ such that $p+2$ is a ${\mathcal P}^+$-integer. Unfortunately, there is no quadratic form which plays a similar role for $\mathcal P$-integers.


\section{An alternative evaluation of $c$}\label{sec:Alternative}

For an alternative approach to evaluating $c$ we will need the $L$-function
\begin{equation}
L(s)=L(s, \chi):=\sum_{n=1}^{\infty}\frac{\chi(n)}{n^s}=\prod_q\left(1-\frac{\chi(q)}{q^s}\right)^{-1}
\end{equation}
corresponding to the non-principal character $\chi$ mod~$(4)$ defined by $\chi(n)=0$ for $n$ even, and $\chi(n)=\pm 1$ for $n\equiv\pm 1$ mod~$(4)$. According to~\cite[Exercise 11.8]{EW} this satisfies
\begin{equation}
\zeta(s)L(s)=\zeta_K(s),
\end{equation}
where $\zeta_K$ is the Dedekind zeta function of the field $K={\mathbb Q}(i)$, defined by
\begin{equation}
\zeta_K(s):=\sum_I\frac{1}{N(I)^s}=\prod_P\left(1-\frac{1}{N(P)^s}\right)^{-1}
\end{equation}
where $I$ and $P$ range over all ideals and prime ideals of ${\mathcal O}_K={\mathbb Z}[i]$, and $N(\;)$ denotes their norm.
(To be more specific, for $K={\mathbb Q}(i)$ the prime $q=2$ ramifies, primes $q\equiv 1$ mod~$(4)$ split into two prime ideals with norm $q$, and primes $q\equiv -1$ mod~$(4)$ are inert, each giving one prime ideal with norm $q^2$, so
\[\zeta_K(s)=\left(1-\frac{1}{2^s}\right)^{-1}\cdot\prod_{q\equiv 1(4)}\left(1-\frac{1}{q^s}\right)^{-2}\cdot\prod_{q\equiv -1(4)}\left(1-\frac{1}{q^{2s}}\right)^{-1};\]
by using
\[1-\frac{1}{q^{2s}}=\left(1-\frac{1}{q^s}\right)\left(1+\frac{1}{q^s}\right)\]
one sees that this is equal to $\zeta(s)L(s)$.)

The Riemann zeta function $\zeta$ has a simple pole at $s=1$, with residue $1$ (see~\cite[Theorem 8.28]{EW}, for example).
By~\cite[Theorem 11.7]{EW} the function $\zeta_K$, like $\zeta$, has a simple pole at $s=1$, but with residue $\rho_K$ given by the class number formula
\begin{equation}
\rho_K=\lim_{s\to 1}(s-1)\zeta_K(s)=\frac{2\pi h}{w\sqrt{|D|}}
\end{equation}
where $D=-4$ is the discriminant of $K$, $h=1$ is its class number and $w=4$ is the number of units in ${\mathcal O}_K$, so that $\rho_K=\pi/4$. (There is also a geometric proof of this, counting lattice points in circles and squares: see~\cite{WikiCNF} for an outline.) 
Since $\zeta$ and $\zeta_K$ have meromorphic continuations in a neighbourhood of $1$, it follows that the function $L=\zeta_K/\zeta$ is meromorphic there, with
\begin{equation}
L(1)=\lim_{s\to 1}\frac{\zeta_K(s)}{\zeta(s)}=\rho_K=\frac{\pi}{4}.
\end{equation}

We can now calculate $\alpha$. Comparing the Euler product expansions for $\zeta$ and $L$ shows that
\[\frac{\left(1-2^{-s}\right)\zeta(s)}{L(s)}=\prod_{q\equiv-1(4)}\frac{1+q^{-s}}{1-q^{-s}}
=\prod_{q\equiv-1(4)}\frac{(1-q^{-2s})}{(1-q^{-s})^2}=\frac{F(s)^2}{F(2s)}.\]
Since the functions $1-2^{-s}$, $L(s)$ and $F(2s)$ are analytic and non-zero at $s=1$,
and $\zeta(s)$ has a simple pole there, it follows that $F(s)^2$ also has a simple pole at $s=1$, so that $\alpha=\frac{1}{2}$. Thus
\begin{equation}
g(x)\sim\frac{cx}{\sqrt{\ln x}}\quad\hbox{as}\quad x\to\infty
\end{equation}
for some constant $c$.
 
We can now determine $c$. In Raikov's Theorem we have
\[c=\frac{H(1)}{\Gamma(1-\alpha)}=\frac{H(1)}{\Gamma(\frac{1}{2})}=\frac{H(1)}{\sqrt\pi},\]
where $H(s)=\sqrt{s-1}F(s)$ near $s=1$, so that
\[H(1)^2=\lim_{s\to 1}\,(s-1)F(s)^2=\frac{F(2)}{2L(1)}\lim_{s\to 1}\,(s-1)\zeta(s)=\frac{2F(2)}{\pi}\]
since $L(1)=\pi/4$ and $\lim_{s\to1}(s-1)\zeta(s)=1$. Thus
\[c=\frac{H(1)}{\sqrt\pi}=\frac{1}{\pi}\sqrt{2F(2)}\]
where
\[F(2)=\prod_{q\equiv -1(4)}\left(1-\frac{1}{q^2}\right)^{-1}=\sum_{n=1}^{\infty}\frac{g_n}{n^2}=1+\frac{1}{3^2}+\frac{1}{7^2}+\cdots.\]

\medskip

The partial sums $\sum_{n=1}^Ng_n/n^2$ of this series converge to $F(2)$ from below as $N\to\infty$.
To obtain a sequence converging to $F(2)$ from above we use
\[\sum_{n=1}^{\infty}\frac{1}{n^2}=\zeta(2)=\frac{\pi^2}{6}\]
(see~\cite[Exercise 9.7]{EW}, for example), so
\[\sum_{{\rm odd}\,n}\frac{1}{n^2}=\zeta(2)-\frac{1}{2^2}\zeta(2)=\frac{\pi^2}{8}\]
and hence
\[F(2)=\frac{\pi^2}{8}-\sum_{n=1}^{\infty}\frac{h_n}{n^2}=\frac{\pi^2}{8}-\left(\frac{1}{5^2}+\frac{1}{13^3}+\cdots\right)\]
where $h_n=1-g_n$ or $0$ as $n$ is odd or even, so that the sum is over all odd $n\not\in{\mathcal P}^*$. Thus
\begin{equation}\label{eq:bounds}
L:= \sum_{n=1}^N\frac{g_n}{n^2}<F(2)<\frac{\pi^2}{8}-\sum_{n=1}^N\frac{h_n}{n^2} \;:= U
\end{equation}
for all $N$, with both bounds converging monotonically to $F(2)$ as $N\to\infty$, so one can evaluate $F(2)$ to any desired accuracy by taking $N$ sufficiently large. 

For example, taking $N=10^9$ in (\ref{eq:bounds}) and evaluating terms to $30$ decimal places gives upper and lower bounds
\[U= 1.16807558580668082351574608710,\]
\[L = 1.16807558530668082351574605022,\]
with $U-L= 5.0000000000000003688 \times 10^{-10}$. For the constant
\[c=\frac{1}{\pi}\sqrt{2F(2)}\]
this gives upper and lower bounds
\[U' = 0.486519888468395460325765388132,\]
\[L' = 0.486519888364266948773906811661,\]
with $U'-L'=1.04128511551858576471 \times 10^{-10}$. The value
\[c=0.4865198883858909971272456405868234055382\ldots\]
given in~(\ref{eq:c_0-value}) lies between these two bounds, and they agree with it in their first nine significant figures.

\medskip

\noindent{\bf Remark} The values of $L$ and $U$ in (\ref{eq:bounds}) for any given $N$ depend on partial sums of the Dirichlet series for the sequences $(g_n)$ and $(h_n)$, representing odd integers $n$ which are or are not $\mathcal P$-integers. For small $N$ the former predominate, so that the lower bounds $L$ are closer to $F(2)$ than the upper bounds $U$ are; since both converge monotonically to $F(2)$, this means that $U$ decreases faster than $L$ increases, especially as non-$\mathcal P$-integers eventually start to predominate as $N$ increases.

To illustrate this, let us denote by $L_k$ and $U_k$ the lower and upper bounds for $F(2)$ computed over the segment $n\le N=10^k$, and let $R_k$ denote the ratio between the values of their rate of change, that is,
\[R_k=\frac{L_k-L_{k-1}}{U_{k-1}-U_k}.\]
For comparison, let $r_k$ denote the ratio $g(x)/h(x)$ of $\mathcal P$-integers to non-$\mathcal P$-integers among the odd integers $n \le x = 10^k$. These ratios $R_k$ and $r_k$ are shown, for $k = 1,\ldots, 9$, in Table~\ref{tab:ratios}.
It follows from Theorem~\ref{th:g-dist} that $R_k$ and $r_k$ converge to $0$ as $k\to\infty$, and Table~\ref{tab:ratios} shows that they do so `in step'.

\begin{table}[htbp]
\begin{center}
\begin{tabular}{c|c|c}
$k$		& $R_k$	& $r_k$ 	 \\
\hline
1      & 3.5966238347        & 4.0000000000               \\
2	& 0.9856093680	& 1.0833333333		\\
3	& 0.6579817478	& 0.6722408026		\\
4	& 0.5213884648	& 0.5060240963		\\
5	& 0.4445379068       & 0.4245825972		\\
6	& 0.3809371214       & 0.3703285499		\\
7	& 0.3382979524        & 0.3309801604		\\
8	& 0.3076691054        & 0.3016477220	         \\
9	& 0.2834484269        & 0.2785807156	         \\
\end{tabular}
\medskip
\caption{Comparison of ratios $R_k$ and $r_k$.}
\label{tab:ratios}
\end{center}
\end{table}


\section{The twin primes conjecture}\label{sec:twin}

It is useful now to recall the heuristic arguments used by Hardy and Littlewood to justify their work in~\cite{HL} on, among various other problems, the Twin Primes Conjecture, which asserts that there are infinitely many pairs of twin primes $p, p+2$.

By the Prime Number Theorem, the probability of a natural number $p\le x$ being prime is asymptotic to $x/\ln x$, so the numbers of primes and of pairs of primes $p, p'\le x$ are asymptotically equivalent to
\[{\rm Li}(x):=\int_2^x\negthickspace\frac{dt}{\ln t}\quad\hbox{and}
\quad I_2(x):=\int_2^x\negthickspace\frac{dt}{(\ln t)^2}.\]
However, if $p$ and $p'$ are related by an equation such as $p'=p+2$, rather than chosen independently, the estimate $I_2(x)$ will be incorrect. Since natural numbers are uniquely determined by their residues mod~$(q)$ for all primes $q$ it makes sense to apply the restriction $p'=p+2$ through its effect on these residues. For each prime $q$, we need to consider the probabilities that $p$ and $p'$, chosen independently or with $p'=p+2$, being both prime, are both coprime to $q$ (we neglect the vanishingly small probability that either of them is equal to $q$). By a simple count of congruence classes mod~$(q)$ we see that these probabilities are respectively
\[\left(1-\frac{1}{q}\right)^2\quad\hbox{and}\quad 1-\frac{\omega_f(q)}{q},\]
where $\omega_f(q)$ is the number of roots of the polynomial $f(t)=t(t+2)$ mod~$(q)$. In order to replace the first probability, implicit in the estimate $I_2(x)$, with the second more correct probability, we therefore multiply $I_2(x)$ by the correction factor
\[C_2(q)=\left(1-\frac{1}{q}\right)^{-2}\left(1-\frac{\omega_f(q)}{q}\right).\]
Doing this for each prime $q$ we multiply $I_2(x)$ by a correction term, called a {\sl Hardy--Littlewood constant},
\[C(f)=\prod_{q\;{\rm prime}}C_2(q)=\prod_{q\;{\rm prime}}\left(1-\frac{1}{q}\right)^{-2}\left(1-\frac{\omega_f(q)}{q}\right).\]

Clearly $\omega_f(2)=1$ while $\omega_f(q)=2$ for each prime $q>2$, so we have an estimate for the number $\pi_2(x)$ of twin prime pairs $p, p+2\le x$ of the form
\begin{equation}\label{eq:E2}
E_2(x)=C(f)I_2(x)=2C_2I_2(x)=2C_2\negthickspace\int_2^x\negthickspace\frac{dt}{(\ln t)^2}\,,
\end{equation}
where the initial factor $2$ is the correction factor $C_2(2)$ corresponding to the prime $q=2$, and
\begin{equation}\label{eq:C2}
C_2:=\prod_{q>2}\left(1-\frac{1}{q}\right)^{-2}\left(1-\frac{2}{q}\right)=\prod_{q>2}\left(1-\frac{1}{(q-1)^2}\right)=0.6601618158\ldots
\end{equation}
with the product over all odd primes $q$.

The conjecture is that
\[\pi_2(x)\sim E_2(x)\quad\hbox{as}\quad x\to\infty.\]
Although it is unproved, there is strong evidence for this conjecture. For example, it is known that
\[\pi_2(10^{18})=808\,675\,888\,577\,436,\]
while Maple evaluates
\[E_2(10^{18})=808\,675\,901\,493\,606.3\ldots.\]
The relative error is $0.0000016\%$.

In 1962 Bateman and Horn~\cite{BH} made a wide-ranging generalisation of the Hardy--Littlewood estimate,
giving a similar conjectured estimate $E(x)$ for the number $Q(x)$ of natural numbers $t\le x$ at which a given finite set of polynomials
$f_1(t),\ldots, f_k(t)\in{\mathbb Z}[t]$
simultaneously take prime values. It assumes that these polynomials satisfy the following conditions, which are obviously necessary for there to be infinitely many such $t$ (Schinzel's Hypothesis~\cite{SS} asserts that they are also sufficient):
\begin{enumerate}
\item each $f_i(t)$ has a positive leading coefficient;
\item each $f_i(t)$ is irreducible in ${\mathbb Z}[t]$;
\item the product $f(t):=f_1(t)\ldots f_k(t)$ is not identically zero modulo any prime.
\end{enumerate}
The Bateman--Horn Conjecture (BHC) asserts that
\begin{equation}\label{eq:BHC}
Q(x)\sim E(x):=C\negthickspace\int_a^x\negthickspace\frac{dt}{\ln f_1(t)\ldots\ln f_k(t)}\quad\hbox{as}\quad x\to\infty,
\end{equation}
where $C$ is the {\sl Hardy--Littlewood constant}
\begin{equation}\label{eq:HL}
C=C(f_1,\ldots, f_k)=\prod_{q\;{\rm prime}}\left(1-\frac{1}{q}\right)^{-k}\left(1-\frac{\omega_f(q)}{q}\right),
\end{equation}
$\omega_f(q)$ is the number of roots of $f$ mod~$(q)$, and
the lower limit $a$ in the integral in (\ref{eq:BHC}) is chosen
so that the integral avoids singularities where $\ln f_i(t)=0$ for some $i$.
If conditions (1) to (3) are satisfied then the infinite product in~(\ref{eq:HL}) converges to a limit $C>0$, while the definite integral in~(\ref{eq:BHC}) diverges to $+\infty$ with $x$, so if the BHC is true then $Q(x)\to+\infty$ and Schinzel's Hypothesis, that the polynomials $f_i$ simultaneously take prime values for infinitely many $t\in\mathbb N$, is verified.

This includes the cases of the twin primes and the Sophie German primes conjectures, where $f_1(t)=t$ and $f_2(t)=t+2$ or $2t+1$ respectively.
The BHC has been proved only in the case of a single polynomial of degree~1:
this is the quantified version, due to de la Vall\'ee Poussin, of Dirichlet's Theorem on primes in an arithmetic progression $at+b$.
Nevertheless, the estimates produced by the BHC, in a wide range of applications, agree remarkably well with observed counts obtained by primality-testing.
As an example, with $f_1(t)=4t+1$ and $f_2(t)=4t+3$, Table~\ref{tab:twin} shows the BHC estimates $E_2^+(x)$ for the numbers $\pi_2^+(x)$ of twin primes $p, p+2\le x$ with $p\equiv 1$ mod~$(4)$, those yielding BCC primes $p$.

\begin{table}[htbp]
\begin{center}
\begin{tabular}{c|c|c|c}
$x$ 			& $\pi_2^+(x)$ & $E_2^+(x)$ & error \\
\hline
10 			& 1 				& 1.148985018 		& 14.8985\,\% \\
$10^2$ 		& 4 				& 5.498634634 		& 37.4659\,\% \\
$10^3$ 		& 19 			& 21.62864106 		& 13.8350\,\% \\
$10^4$ 		& 105 			& 105.8363607 		& 0.79653\,\% \\
$10^5$ 		& 604 			& 623.0852586 		& 3.15981\,\% \\
$10^6$ 		& 4\,046 		& 4\,122.745734 		& 1.89683\,\% \\
$10^7$ 		& 29\,482 		& 29\,375.63915 		& $-0.3608\,\%$ \\
$10^8$ 		& 220\,419 		& 220\,182.6280 		& $-0.1072\,\%$ \\
$10^9$ 		& 1\,712\,731 	& 1\,712\,652.809 	& $-0.0046\,\%$ \\
$10^{10}$ 	& 13\,706\,592	& 13\,705\,706.99	& $-0.0065\,\%$ \\
\end{tabular}
\medskip
\caption{Numbers of pairs of twin primes $p,p+2\le x$ with $p\equiv 1$ mod~$(4)$.}
\label{tab:twin}
\end{center}''
\end{table}


\section{In the footsteps of Hardy and Littlewood}\label{sec:footsteps}

We will now adapt the heuristic arguments used by Hardy and Littlewood, and later by Bateman and Horn, to the slightly different context of our BCC primes problem.

Using the Prime Number Theorem and Theorem~\ref{th:g-dist} to give the distributions of prime numbers $p$ and of $\mathcal P$-integers $n$, we first consider
\[I(x):=\int_2^x\negthickspace\frac{c\,dt}{(\ln t)^{3/2}}\]
as an estimate for the number $a(x)$ of BCC primes $p\le x$, where $c$ is as in (\ref{eq:c_0-value}). Of course, this treats $p$ and $n$ as independent random variables, and takes no account of the fact that $n=p+2$. In following the example of Hardy and Littlewood, with a similar heuristic justification, we now apply correction factors only for primes $q\not\in{\mathcal P}$, that is, $q=2$ or $q\equiv 1$ mod~$(4)$, since it is for such $q$ that we need to replace the probability that independent variables $p$ and $n$ are both coprime to $q$ with the corresponding probability for $p$ and $p+2$. We therefore multiply this estimate $I(x)$ by $2C_2^+$ where
\begin{equation}\label{eq:C2+}
C_2^+:=\prod_{q\equiv 1(4)}\left(1-\frac{1}{q}\right)^{-2}\left(1-\frac{2}{q}\right)
=\prod_{q\equiv 1(4)}\left(1-\frac{1}{(q-1)^2}\right),
\end{equation}
with the extra factor $2$ corresponding to the prime $q=2$.

This leads to the following conjecture:

\begin{conj}\label{conj:a-distr}
The function $a(x)$ which counts BCC primes $p\le x$ satisfies
\begin{equation}\label{eq:a-distr}
a(x)\sim 2C_2^+I(x)
=c'\negthickspace\int_2^x\negthickspace\frac{dt}{(\ln t)^{3/2}}\quad\hbox{as}\quad x\to\infty,
\end{equation}
where
\begin{equation}\label{eq:c-const}
c'=2C_2^+c=\prod_{q\equiv 1\,(4)}\left(1-\frac{1}{(q-1)^2}\right)\cdot
\sqrt{\frac{e^{\gamma}}{\pi}}\cdot C(4,1).
\end{equation}
\end{conj}


\section{Computations}\label{sec:Computations}

The various estimates we have discussed can be tested computationally by using Maple. The definite integrals, such as those appearing in (\ref{eq:E2}) and (\ref{eq:a-distr}), are evaluated accurately and rapidly, even for large values of $x$, by numerical integration. The Hardy--Littlewood constants, such as those in~(\ref{eq:C2}) and (\ref{eq:C2+}), are defined as infinite products which converge conditionally and slowly; good approximations can be found by taking partial products over large initial segments of the relevant primes, in increasing order, but on a laptop these calculations can take a matter of hours.
Having evaluated the various estimates, one can compare them with the actual numbers of terms being counted by using the Rabin--Miller primality test within Maple.

Using Maple to evaluate (\ref{eq:C2+}) for primes $q\le 10^8$ we find that
\[C_2^+=0.9230611326\ldots\]
The values of $\sqrt{e^{\gamma}/\pi}$ and $C(4,1)$ given in Section~\ref{sec:Distribution} allow us to deduce that
\begin{equation}\label{eq:c-val}
c'=0.8981751984\ldots
\end{equation}

To test the value for $c'$ in~(\ref{eq:c-val}) let us define $C(x)$ by
\[a(x)=C(x)\negthickspace\int_2^x\negthickspace(\ln t)^{-3/2}dt,\]
so we expect that $C(x)\to c'$ as $x\to\infty$.
Using Maple to evaluate $a(x)$ and $\int_2^x(\ln t)^{-3/2}dt$, we found the values of $C(x)$ shown in Tables~\ref{tab:constants1} and~\ref{tab:constants2}. We note that, after some initial instability due to the relatively small numbers of BCC primes appearing, the values of $C(x)$ decrease almost monotonically towards $c'$ from $x=10^7$ to $x=35\cdot 10^9$.


\begin{table}[htbp]
\begin{center}
\begin{tabular}{c|c|c}
$x$				& $a(x)$				& $C(x)$		\\
\hline
$10^1$			& 2					& 0.4688555840	\\
$10^2$			& 12					& 0.7153107401	\\
$10^3$			& 65					& 0.8472363117	\\
$10^4$			& 388				& 0.8706477077	\\
$10^5$			& 2\,708				& 0.9004062930	\\
$10^6$			& 19\,969			& 0.9024565742  	\\
$10^7$			& 155\,369			& 0.9040719795	\\
$10^8$			& 1\,250\,182		& 0.9023589943	\\
$10^9$			& 10\,345\,920		& 0.9011815839	\\
$10^{10}$		& 87\,545\,946		& 0.9010054301	\\
\end{tabular}
\medskip
\caption{Evolution of the coefficients $C(x)$ for $x=10^k$, $k=1,\ldots,10$.}
\label{tab:constants1}
\end{center}
\end{table}


\begin{table}[htbp]
\begin{center}
\begin{tabular}{c|c|c}
$x$				& $a(x)$				& $C(x)$		\\
\hline
$11\cdot 10^9$	& 95\,675\,252		& 0.9010083668	\\
$12\cdot 10^9$	& 103\,758\,501		& 0.9010266288	\\
$13\cdot 10^9$	& 111\,794\,166		& 0.9010090565	\\
$14\cdot 10^9$	& 119\,795\,477		& 0.9010345605	\\
$15\cdot 10^9$	& 127\,757\,388		& 0.9010373485	\\
$16\cdot 10^9$	& 135\,683\,004		& 0.9010241227	\\
$17\cdot 10^9$	& 143\,578\,133		& 0.9010192899	\\
$18\cdot 10^9$	& 151\,444\,525		& 0.9010197264	\\
$19\cdot 10^9$	& 159\,283\,669		& 0.9010228308	\\
$20\cdot 10^9$	& 167\,092\,278		& 0.9010018195	\\
$21\cdot 10^9$	& 174\,878\,590		& 0.9009949342	\\
$22\cdot 10^9$	& 182\,641\,563		& 0.9009888293	\\
$23\cdot 10^9$	& 190\,382\,424		& 0.9009836658	\\
$24\cdot 10^9$	& 198\,100\,186		& 0.9009699991       \\
$25\cdot 10^9$	& 205\,796\,174		& 0.9009503295	\\
$26\cdot 10^9$	& 213\,474\,480		& 0.9009388161	\\
$27\cdot 10^9$	& 221\,136\,287		& 0.9009360946	\\
$28\cdot 10^9$	& 228\,778\,913		& 0.9009277506	\\
$29\cdot 10^9$	& 236\,403\,538		& 0.9009161678	\\
$30\cdot 10^9$	& 244\,014\,302		& 0.9009145456	\\
$31\cdot 10^9$	& 251\,607\,900		& 0.9009079931	\\
$32\cdot 10^9$	& 259\,184\,139		& 0.9008943241	\\
$33\cdot 10^9$	& 266\,743\,992		& 0.9008756612	\\
$34\cdot 10^9$	& 274\,293\,778		& 0.9008715586	\\
$35\cdot 10^9$	& 281\,827\,468		& 0.9008601041
\end{tabular}
\medskip
\caption{Evolution of the coefficients $C(x)$ for $x=m\cdot 10^9$, $m=11,\ldots,35$.}\label{tab:constants2}
\end{center}
\end{table}

Beyond this point, in view of the modest computing facilities available we considered segments $[x,y]\subset\mathbb R$ of length $y-x=10^8$ or $10^7$, starting at values $x=10^{15}, 10^{20},\ldots, 10^{50}$. The results are shown in Table~\ref{tab:segments}, where $a(x,y)=a(y)-a(x)$ is the number of BCC primes $p\in [x,y]$, and $C(x,y)$ is defined by the equation
\[a(x, y) = C(x,y) \cdot\negthinspace \int_x^y\negthickspace \frac{dt}{(\ln t)^{3/2}}.\]

Initially we see a further monotonic decrease towards $c'$, but then $C(x,y)<c'$ when $x=10^{30}$ and $y=10^{30}+10^8$. We are not very concerned about this: for instance, Littlewood~\cite{Lit} famously showed that the error ${\rm Li}(x)-\pi(x)$ in the Prime Number Theorem changes sign infinitely many times (see also~\cite{Rob} for a similar phenomenon related to Mertens's Third Theorem), so why not here? The instability in the last four rows is perhaps due to the relatively small numbers of BCC primes appearing in these shorter intervals: for instance, compare the values of $a(x,y)$ with that of $a(x)$ for the equivalent interval $[0,10^7]$.

\begin{table}[htbp]
\begin{center}
\begin{tabular}{c|c|c|c}
$x$			& $y-x$	 & $a(x,y)$	& $C(x,y)$				\\
\hline
$10^{15}$	& $10^8$	 & 442\,649		& 0.8985037831	\\
$10^{20}$	& $10^8$	 & 287\,429		& 0.8982538979	\\
$10^{25}$	& $10^8$	 & 205\,949		& 0.8994835673	\\
$10^{30}$	& $10^8$	 & 156\,398		& 0.8979178604	\\
$10^{35}$	& $10^7$	 & 12\,389		& 0.8963174519	\\
$10^{40}$	& $10^7$	 & 10\,299		& 0.9103503807	\\
$10^{45}$	& $10^7$	 & 8\,554		& 0.9022181293	\\
$10^{50}$	& $10^7$	 & 7\,507		& 0.9273527313	\\
\end{tabular}
\medskip
\caption{Estimation over segments $[x,y]$.}\label{tab:segments}
\end{center}
\end{table}


\section{Some remarks}\label{sec:remarks}

{\bf 1.} The estimate $g(x)\sim cx/(\ln x)^{1/2}$ in Theorem~\ref{th:g-dist}, including the evaluation of $c$ in~(\ref{eq:c_0-value}), is proved, as a consequence of Raikov's Theorem. The estimate $a(x)\sim c'\int_2^x(\ln t)^{-3/2}dt$ in Conjecture~\ref{conj:a-distr}, including the formula for $c'$ given there, is just a conjecture, and we do not expect to see a proof for it soon. Nevertheless, we feel that the computational data in Section~\ref{sec:Computations} give plausible evidence for the validity of this estimate, and hence for the infinitude of BCC primes.

\smallskip

{\bf 2.} In the construction of Riemann surfaces used in~\cite{BCC} and considered here, it is important for group-theoretic reasons that the prime $p$ and the $\mathcal P$-integer $n$ should differ by $2$. However, it is clear that in the arguments we have used one could generalise this relationship, and still obtain similar estimates for the distribution of such primes. For instance, one could replace the difference $2$ here with any non-zero even integer. The following section gives a less trivial example, also arising from~\cite{BCC}.


\section{A similar problem}\label{sec:Similar}

In~\cite{BCC}, having obtained an Accola--Maclachlan-type lower bound $M(g)\ge 4(g-1)$ for the largest possible order $M(g)$ of the automorphism group of a pseudo-real Riemann surface of odd genus $g\ge 3$, Bujalance, Cirre and Conder have presented a set $\mathcal G$ of genera $g$ for which this bound is sharp. (The sharpness of the corresponding bound $M(g)\ge 2g$ for even genera $g\ge 2$ is left as an open problem.) For this they give a similar construction (see the remark following their Theorem~5.2) of a family of pseudo-real Riemann surfaces of genus $g=2p + 1$ where $p$ is
what we will call a BCC2 prime, defined as follows:
\begin{defi}
A BCC2 prime is a prime number $p$ such that
\begin{equation}\label{eq:congs}
p\equiv 3\; {\rm mod}\,(8),\; p\equiv 2\;\hbox{or}\; 5\; {\rm mod}\,(9)\;\hbox{and}\; p\not\equiv 5\;{\rm mod}\,(7),
\end{equation}
and
\begin{equation}\label{eq:divisors}
n:=p+1\;\hbox{is not divisible by}\;11, 23, 47\;\hbox{or any prime} \; q\equiv 1\;\hbox{mod}\;(3)\;\hbox{or}\; q\equiv 1\;\hbox{mod}\;(4).
\end{equation}
\end{defi}

This raises the following problem, which appears to be even more challenging than Problem~\ref{pro:BCC}:
\begin{prob}\label{pro:BCC2}
Are there infinitely many BCC2 primes?
\end{prob}

These primes are much rarer than the BCC primes considered earlier:
for example, compare Table~\ref{tab:BCC2primes}, which shows the first sixteen of them, with the corresponding list of BCC primes in~(\ref{eq:smallBCC}).

\begin{table}[htbp]
\begin{center}
\begin{tabular}{c|c|c|c}
$p$		& factors of $n=p+1=12m$  & $m$ mod~$(7)$	& $r=m$ mod~$(84)$  \\
\hline
$11$	& $2^2\cdot 3$	 & 		& 	\\
$1283$	& $2^2\cdot 3\cdot 107$	 & 	$2$	& $23$	\\
$1571$	& $2^2\cdot 3\cdot 131$	 & 	$5$	& $47$	\\
$2003$	& $2^2\cdot 3\cdot 167$	 & 	$6$	& $83$	\\
$3011$	& $2^2\cdot 3\cdot 251$	 & 	$6$	& $83$	\\
$7043$	& $2^2\cdot 3\cdot 587$	 & 	$6$	& $83$	\\
$7907$	& $2^2\cdot 3\cdot 659$	 & 	$1$	& $71$	\\
$8627$	& $2^2\cdot 3\cdot 719$	 & 	$5$	& $47$	\\
$9923$	& $2^2\cdot 3\cdot 827$	 & 	$1$	& $71$	\\
$10\,067$	& $2^2\cdot 3\cdot 839$	 & 	$6$	& $83$	\\
$15\,107$	& $2^2\cdot 3\cdot 1259$	 & 	$6$	& $83$	\\
$15\,683$	& $2^2\cdot 3\cdot 1307$	 & 	$6$	& $83$	\\
$17\,123$	& $2^2\cdot 3\cdot 1427$	 & 	$6$	& $83$	\\
$17\,987$	& $2^2\cdot 3\cdot 1499$	 & 	$1$	& $71$	\\
$18\,131$	& $2^2\cdot 3\cdot 1511$	 & 	$6$	& $83$	\\
$19\,427$	& $2^2\cdot 3\cdot 1619$	 & 	$2$	& $23$	\\
\end{tabular}
\medskip
\caption{The first sixteen BCC2 primes $p$.}
\label{tab:BCC2primes}
\end{center}
\end{table}

First let us restate the definition of BCC2 primes in terms of $n$. The congruences mod~$(8)$ and mod~$(9)$ in condition~(\ref{eq:congs}) are equivalent to $n=p+1\equiv 4$ mod~$(8)$ and $n\equiv 3$ or $6$ mod~$(9)$, that is, $n=12m$ for some $m$ coprime to $12$; when this is satisfied the condition $p\not\equiv 5$ mod~$(7)$ in~(\ref{eq:congs}) is equivalent to $m\not\equiv 4$ mod~$(7)$. Thus condition~(\ref{eq:congs}) is equivalent to
\begin{equation}\label{eq:newcongs}
n=12m\;\hbox{where}\;(12, m)=1\;\hbox{and}\; m\not\equiv 4\,{\rm mod}\,(7).
\end{equation}
Similarly, when this is satisfied, condition~(\ref{eq:divisors}) is equivalent to
\begin{equation}\label{eq:newdivisors}
\hbox{each prime $q$ dividing $m$ satisfies}\; 47<q\equiv -1\, {\rm mod}~(12).
\end{equation}
Thus conditions~(\ref{eq:congs}) and~(\ref{eq:divisors}), taken together, are equivalent to (\ref{eq:newcongs}) and (\ref{eq:newdivisors}), also taken together. Recall that in addition we require $n-1$ ($=p$) to be prime (and greater than $7$), which implies that $m\not\equiv 3$ mod~$(7)$ since otherwise
we would have $12m-1\equiv 36-1 \equiv 0$ mod $(7)$, while conditions (\ref{eq:newcongs}) and (\ref{eq:newdivisors}) imply that $m\not\equiv 4$ or $0$ mod~$(7)$, so we have
\begin{equation}\label{eq:mmod7}
m\equiv 1, 2, 5\;\hbox{or}\;6\;\hbox{mod}\;(7).
\end{equation}

In Table~\ref{tab:BCC2primes}, apart from $p=11$, all entries have the property that $n=12m$ for some prime $m$. Indeed, the first exception to this is $p=41\,771$, with $n=12\cdot 59^2$. Moreover, among the $172\,515$ BCC2 primes $p\le 10^9$, there are only $41\,711$ for which $m$ is not prime. Motivated by this, 
instead of attempting to estimate the distribution of the whole set of BCC2 primes, as we did in the case of the BCC primes, we will apply the Bateman--Horn Conjecture (see Section~\ref{sec:twin}) to give strong evidence that there are infinitely many of them for which $m$ is prime. (This is analogous to the argument at the end of Section~\ref{sec:twin}, where we considered those BCC primes $p$ which are members of twin primes $p, p+2$.)

Let us therefore assume that $n=12m$ for some prime $m>47$ such that $m\equiv -1$ mod~$(12)$ and $m\equiv 1, 2, 5$ or $6$ mod~$(7)$, or equivalently $m\equiv r$ mod~$(84)$ where $r=71, 23, 47$ or $83$ by the Chinese Remainder Theorem. Clearly, any such choice of $m$ satisfies conditions (\ref{eq:newcongs}), (\ref{eq:newdivisors}) and (\ref{eq:mmod7}). For each such $r$ we are therefore looking for integers $t\ge 0$ such that the polynomials \[f_1(t)=84t+r\]
and
\[f_2(t)=12f_1(t)-1=1008t+12r-1\]
both take prime values, namely $m$ and $p$. These two polynomials are irreducible, with positive leading coefficients, and their product $f$ is not identically zero modulo any prime, so we have four instances of the Bateman--Horn Conjecture, one for each value of $r$. In each case the BHC gives an asymptotic estimate $E_r(x)$ for the number $Q_r(x)$ of integers $t\le x$ such that $f_1(t)$ and $f_2(t)$ are both prime. We have $\omega_f(q)=0$ for the primes $q=2, 3$ and $7$, and $\omega_f(q)=2$ for all other primes, so the Hardy--Littlewood constant $C=C(f_1,f_2)$, which is independent of $r$, is positive. In fact, comparing the infinite products for $C$ and for $C_{\rm twins}=2C_2$ (see Section~\ref{sec:twin}), which differ at only these three primes, shows that
\[ C=\frac{42}{5}\cdot C_{\rm twins}=\frac{42}{5}\cdot 1.3203236316\ldots=11.0907185062\ldots .\]
It follows that $E_r(x)\to+\infty$ as $x\to\infty$, giving evidence that for each $r$ there are infinitely many prime values $m$ and $p$ for $f_1$ and $f_2$, and thus infinitely many BCC2 primes $p$.

To test this approach we found that of the $172\,515$ BCC2 primes $p\le 10^9$, the number with $m$ prime, as we have been assuming, is $Q=130\,804$.
The four BHC estimates
\[E_r(x)=C\cdot\negthinspace\int_2^x\negthickspace\frac{dt}{\ln f_1(t)\ln f_2(t)} \]
were computed by taking $x=992\,062$; beyond this point the primes $p=1008t+12r-1$ begin to be greater than $10^9$. The four values obtained were almost identical, ranging from $32\,709.13$ to $32\,709.29$. Their sum gave an estimate of $E=130\,836.81$ for $Q$, with an error of $0.025\%$.


\section{Appendix: the pseudo-real surfaces}\label{sec:pseudo-real}

Here we briefly summarise the construction of the pseudo-real Riemann surfaces in~\cite{BCC} which give rise to Problems~\ref{pro:BCC} and \ref{pro:BCC2}. In the first case, for any even $g\ge 2$ let $\Gamma$ be an NEC group with signature $(1;-;[2,2,g];\{-\})$. This means that $\Gamma$ acts as a group of isometries of the hyperbolic plane $\mathbb H$, with canonical generators $d$ (a glide reflection) and $x_1, x_2, x_3$ (elliptic elements), and defining relations
\[x_1^2=x_2^2=x_3^g=d^2x_1x_2x_3=1,\]
so that the quotient-surface ${\mathbb H}/\Gamma$ is the real projective plane with three cone-points of orders $2,2$ and $g$.

Now let $\theta:\Gamma\to {\rm C}_{2g}=\langle u\mid u^{2g}=1\rangle$ be the surface-kernel epimorphism defined by
\[d\mapsto u,\; x_1\mapsto u^g,\; x_2\mapsto u^g,\; x_3\mapsto u^{-2},\]
and let $K=\ker\,\theta$. Then $K$ is contained in the orientation-preserving subgroup $\Gamma^+$ of index $2$ in $\Gamma$, so $S:={\mathbb H}/K$ is a compact Riemann surface, of genus $g$ by the Riemann--Hurwitz formula. There is a natural action of $ \Gamma/K$ and hence of ${\rm C}_{2g}$ on $S$, with the elements of the subgroup $\langle u^2\rangle\cong {\rm C}_g$ and its complement acting as conformal and anticonformal automorphisms of $S$, and with $S/{\rm C}_{2g}\cong{\mathbb H}/\Gamma$; the canonical double cover $S/{\rm C}_g\cong{\mathbb H}/\Gamma^+$ of this orbifold is the Riemann sphere with four cone-points of order $2$ and two of order $g$.

The generator $u$ of ${\rm C}_{2g}$ acts as an anticonformal automorphism of $S$ of order $2g$, whereas the only involution in $G$ (namely $u^g$) acts conformally. One can choose the cone-points so that $\Gamma$ is a maximal NEC group, in which case $S$ has no further automorphisms (conformal or anticonformal), and is therefore pseudo-real. A detailed argument in~\cite[Theorem~4.8]{BCC} shows that provided $g$ satisfies various other conditions (namely those defining the set $\mathcal A$, together with $g>30$), no pseudo-real surface of genus $g$ can have an abelian group of conformal automorphisms with more than $g$ elements, so that $M^+_{\rm ab}(g)=g$.

The construction of the surfaces $S$ giving rise to Problem~\ref{pro:BCC2} is similar, except that in this case $\Gamma$ has signature $(1;-;[2,2,2];\{-\})$, $S$ has genus $g=2p+1$ with a full automorphism group
\[\langle u, v \mid u^2=v^{4p}=1, v^u=v^{2p-1}\rangle\cong{\rm C}_{4p}\rtimes{\rm C}_2\]
of order $8p=4(g-1)$, and $\theta$ is given by
\[d\mapsto uv,\; x_1\mapsto u,\; x_2\mapsto u,\; x_3\mapsto v^{2p}.\]
It is shown in~\cite[Theorem~5.3]{BCC} that if $p$ is a prime satisfying the conditions stated in Section~\ref{sec:Similar} then no pseudo-real surface of genus $g$ has a larger automorphism group, so that $M(g)=4(g-1)$.


\medskip

\noindent{\bf Acknowledgments} The authors are grateful to Emilio Bujalance, Javier Cirre and Marston Conder for creating the problems addressed in this paper, and to Ken Ford for some very helpful advice.


\end{document}